\documentclass{amsart}
\usepackage[all]{xy}
\usepackage[T1]{fontenc}
\usepackage{times}
\usepackage{amssymb,epsfig,verbatim,xy}
\usepackage{enumitem}
\theoremstyle{plain}

\usepackage[T1]{fontenc}
\usepackage[utf8]{inputenc}
\usepackage[english]{babel}
\usepackage{amssymb}
\usepackage{amsthm}
\usepackage{graphics}
\usepackage{amsmath}
\usepackage{amstext}
\usepackage{multirow}
\usepackage{hyperref}

\usepackage{graphicx}
\usepackage{amsmath}
\usepackage{amssymb}
\usepackage{amsthm}
\usepackage{amstext}
\usepackage{epstopdf}
\usepackage{mathrsfs}
\usepackage{amscd}
\usepackage{color}

\theoremstyle{remark}

\title {A non-algebraizable foliation}
\date{\today}

\author{Paulo Sad }

\begin {document}

\maketitle

\begin{abstract}
It is presented an example of a holomorphic  foliation of a non-algebraizable surface which is topologically equivalent to an algebraic foliation. 

\end{abstract}

\vspace {0.3 in}
In this note we show in a specific example how to change an algebraic foliation to another holomorphic foliation which is topologically but not holomorphically equivalent. Of course several examples of moduli spaces of foliations are known, the simplest one appears with local 2-dimensional singularities in the Poincaré domain; there is also an extensive analysis in \cite{Ma} for more general singularities. The difference here is that we change the holomorphic type of the surface where the foliation is defined, but preserving the existence of a similar (from the topological point of view) holomorphic foliation of the new surface. More precisely, we exhibit  a pair of holomorphic foliations $\mathcal F$ and $\mathcal G$ with the following properties:

\vspace {0.1 in}
\noindent (i) $\mathcal F$ is defined in the projective plane ${\mathbb C}P^2$; it has an invariant smooth compact curve $C$ of genus 3 with 16 dicritical singularities. The local holonomy map relatively to $C$ at each singularity is the identity, but the holonomy group of $C$ has infinitely many elements.%

\vspace {0.1 in}
\noindent (ii) $\mathcal G$ is defined in some non algebraizable surface containing a copy of $C$ and it is topologically equivalent to the restriction of $\mathcal F$ to some neighborhood of $C$.

\vspace {0.1in} In the language of \cite{Ma}, the moduli space of $\mathcal F$ is not trivial.
  
\vspace {0.1 in}
We use the construction  presented in \cite{F} which allows to go from an algebraic surface to a non algebraizable one, but in a different way since we deal with extending holomorphic applications between (non smooth) Stein spaces. We apply the construction to the surface associated to the Neeman's example which  has a foliation with a genus 3 compact curve as a leaf (\cite {}); that foliation can be transported to the non-algebraizable surface.

\section {The Surfaces}

We start with a smooth plane quartic $C \subset {\mathbb C}P^2$ and select 16 points: two pairs of points belong to two bitangent lines and the remaining 12 points to the intersection with a smooth cubic. Then we blow-up ${\mathbb C}P^2$ at these poins to obtain a surface $\Hat {{\mathbb C}P^2}$; the remarkable fact here is the existence of a smooth strictly psh function $f$ defined on a neighborhood $V =\underset {t> t_0}\bigcup f^{-1}(t) \cup {\Hat C} \subset {\Hat {{\mathbb C}P^2}}$ of the strict transform $\hat C$ of $C$ which tends to $\infty$ as we approach $\hat C$. As a consequence, if $\bar t$ is such that the real 3-manifold $f^{-1}({\bar t})$ is smooth then by \cite{Na} the set ${\Hat {{\mathbb C}P^2}} \setminus \underset{t \geq {\bar t}}\bigcup  f^{-1}(t)$ is holomorphically convex and contains a maximal exceptional subset that can be contracted to a finite number of points \cite {Gra}. We get a (singular) Stein surface $M$ with a finite number of singularities; we denote by $M_{\infty}$ the compact surface $M\cup V$.

\vspace {0.1in}
\noindent {\bf Remark}  $\Hat {{\mathbb C}P^2}\setminus {\Hat C}$ is an example of a 1-convex complex surface,  therefore it is a holomorphically convex set (not necessarily Stein, differently from the construction in Section 6 of \cite{F}).

\vspace {0.1 in}
 \noindent Let us introduce the following notation:

\vspace {0.1 in}
\noindent (i)  $ p_1, \dots,p_{16}$ the points of $C$ where ${\mathbb C}P^2$ is blown-up.

\noindent (ii) $D_1,\dots, D_{16}$ the exceptional divisors.

\noindent (iii) $ \{{\hat p}_j\}= {\Hat C} \cap D_j$

\noindent (iv) $V_j \subset \Hat{ {\mathbb C}P^2}$ neighborhood of $D_j$   that is blown-down to a small neighborhood of $p_j$

\vspace {0.1 in}
\noindent The surface $\Hat S:= V \underset {id}\cup (V_1\cup \dots \cup V_{16}) \subset M_{\infty}$ may be seen as $V$ glued to each $V_j$ by the identity map near $\hat p_j$.

Now we replace one of the identifications, say near $\hat p_1$, by another biholomorphism $\phi$ that sends a smooth curve $l_1\subset V$ transverse to $\Hat C$ and $D_1$ at $\hat {p_1}$ into $D_1$. Let us explain briefly how this glueing is done.

\vspace {0.1 in}
 We take  a bidisk $\{(x,y);|x| \leq r, |y| \leq r\}$ such that $\hat{p_1}=(0,0)$, $\{(x,y);|x| \leq r, y=0\}$ is contained in $\Hat C$  and  $\{(x,y);|x|\leq r, |y|=r\}$ is outside  $\bar V$. We assume that $l_1$ is parametrized as $x=x(y)$ with $|x(y)|< r/2$  and $x(0)=0$. In $V_1$ we take the bidisc  $\{(x^{\prime},y^{\prime}); |x^{\prime}|\leq r, |y^{\prime}|\leq r\}$ such that $\hat {p_1}=(0,0)$, $\{(x^{\prime},y^{\prime}); |x^{\prime}|=r, |y^{\prime} \leq r\}$  is contained in $\partial$$V_1$, $\{(x^{\prime},0); |x^{\prime}| \leq r \} \subset {\Hat C}$ and $\{(0,y^{\prime}); |y^{\prime}|\leq r\} \subset D_1$. The mapping $(x^{\prime},y^{\prime})=\phi (x,y)$ satisfies: $\phi$ sends $|x|=r$ to $|x^{\prime}|=r$, $\phi (x,0)=(x,0)$ , it preserves the horizontals $y^{\prime}=c$  and $\phi(x(c),c)=(0,c)$ for each $|c|\leq r$ ($\phi$ can be taken as a homography for each horizontal). It is easy to see that we can define a continuous isotopy from $\phi| _{\{(x,y) ;|x|=r\}}$ to the identity restricted to $\{(x,y); |x|=r_1\}$  ($r_1$ close to $ r$)  along the horizontals $\{(x,y); r_1\leq |x| \leq  r, |y|=c\}$, $c$ sufficiently small.

\vspace {0,1 in}
We put $\Hat S^{\prime}= (V \underset {\phi}\cup V_1) \underset {id}\cup (V_2
\cup \dots \cup V_{16})$; it is a smooth complex surface. Let us denote by $V^{\prime}$ and $ D_1^{\prime}$ the copies of $V$ and $ D_1$ contained in ${\Hat S}^{\prime}$ and by $i$ the natural biholomorphism which takes $V$ to $V^{\prime}$.

\vspace {0.1 in}
\noindent {\bf Lemma}  Suppose that $ l_1$ is contained in some irreducible algebraic curve $ L_1$  of $\Hat  {{\mathbb C}P^2}$ which intersects $\Hat C$ in another point $\hat p$ different from $\hat {p_1}$. Then there is no embedding of ${\Hat S}^{\prime}$ into a compact holomorphic surface.
\begin {proof} Without any loss of generality we will prove that ${\Hat S}^{\prime}$ is not an open subset of a compact holomorphic surface $T$.
The function $f^{\prime}= f \circ {i}^{-1}$ is also a strictly smooth psh function of $V^{\prime}$. We have $V^{\prime}=\underset {t>t_0}\cup {f^{\prime}}^{-1}(\bar t)$ and ${f^{\prime}}^{-1}(\bar t)$ is a smooth real 3-manifold. Consequently $T\setminus \underset {t\ge {\bar t}}\cup {f^{\prime}}^{-1}(t)$ is a holomorphically convex set with boundary ${f^{\prime}}^{-1}(\bar t)$. We contract the exceptional subset of  $T\setminus \underset{t\ge  {\bar t}}\cup {f^{\prime}}^{-1}(t)$ to get a (singular) Stein surface $M^{\prime}$. We write $M_{\infty}^{\prime} = M^{\prime} \cup V^{\prime}$.
Consider now $i$ restricted to $\underset {t_0<t<{\bar t}}\bigcup f^{-1}(t) \subset M$ and taking values in $\underset {t_0<t<{\bar t}}\bigcup { f^{\prime}}^{-1}(t) \subset M^{\prime}$. By a generalized Hartogs Theorem for holomorphic applications between Stein surfaces (\cite{K}) we can  extend $i$ to a holomorphic application $I$ defined in $M$ and taking values in $M^{\prime}$, and {\it a fortiori} from $M_{\infty}$ to $M^{\prime}_{\infty}$. 
Now we look to $I^{-1}(D_1^{\prime })$ (we keep using the same notation for $D_1^{\prime}$ after the contractions  done in $T\setminus \underset {t\ge {\bar t}}\cup {f^{\prime}}^{-1}(t))$. It contains $L_1$, so that $\{\hat p_1, \hat p\}\subset I^{-1}(D_1^{\prime})$, contradiction.
\end {proof}

As a consequence, $\Hat S$ is not biholomorphically equivalent to ${\Hat S}^{\prime}$ when we use special curves for $l_1$. For the sake of clarity, let us indicate these surfaces as  ${\Hat S_1}^{\prime}$

\vspace {0.1 in}
Let $S_1^{\prime}$ be the surface obtained from ${\Hat S_1}^{\prime}$ blowing-down all the (-1) divisors. The Lemma implies

\vspace {0.1 in}
\noindent {\bf Theorem} The germ of $S_1^{\prime}$ along $C$ is not algebraizable.

\section {The Foliations}

It is known (\cite {S}) that we can select  cubics in the Neeman's construction such that when the blow up's occur at  the 12 points of intersection with the quartic $C$ (besides the points lying in the bitangents) the resulting surface $\Hat {{\mathbb C}P^2}$ has a global foliation  $\Hat {\mathcal F}$
which contains $\hat C$ as a leaf (withouth singularities!). The surface $\Hat S$ we considered in the last Section has therefore a regular foliation $\Hat {\mathcal F}$.

Let us go back to the construction of the surface $\Hat {S^{\prime}}$ of the last Section. In the description of the glueing map $\phi$ we used bidiscs in $V$ and $V_1$; their  horizontals can be supposed to be contained in the leaves of $\Hat {\mathcal F}$. Since $\phi$ preserves horizontals, we have a new foliation $\Hat{ {\mathcal F}^{\prime}}$ defined in $S^{\prime}$. We write ${\mathcal F}_1^{\prime}$ for such a foliation in $S_1^{\prime}$.

\vspace {0.1in}
\noindent {\bf Proposition} $\Hat{ \mathcal F}$ and $\Hat {\mathcal F}^{\prime}$ are topologically equivalent.
\begin {proof} We consider $\Hat S$ as the glueing of $V$ to $V_1$ by the identity. We define a topological equivalence $ \Psi$ as the Identity outside $V$; it follows that in $V_1$ it becomes $\phi$ for $|x|=r$.  We extend the homeomorphism between $|x|=r_1$ and $|x|=r$ as the isotopy we constructed in the last Section and as the Identity for $|x\leq r_1$. Of course we may complete the definition as the Identity for a neighborhood of $D_1$ inside $V_1$.
\end {proof}

\vspace {0.1in}
We remark that the foliations are holomorphically equivalent around the corresponding singularties and  that $\Hat C$ has holomorphically conjugated groups of holonomy.

\vspace {0.2in}
\noindent  {\bf Corollary} The foliations $\Hat {\mathcal F}$ and $\Hat {\mathcal F}_1$ are topologically but not holomorphically equivalent.

\vspace {0.2 in}
\noindent After blowing down the (-1) divisors, we see that $\mathcal F$ and ${\mathcal F}^{\prime}$ are topologically but not holomorphically equivalent.


\vskip 0.1in
\noindent\author{Paulo Sad}
\address{IMPA, Estrada Dona Castorina, 110, Horto, 22460-320,  Rio de Janeiro, Brasil.}

\noindent \email{sad@impa.br}


\end{document}